\documentclass[11pt,a4paper]{amsart}
\usepackage[T1]{fontenc}
\usepackage{lmodern,amsmath,amssymb,amsthm}
\usepackage[margin=26mm]{geometry}
\usepackage{tikz-cd}
\usepackage[hidelinks,pdfusetitle]{hyperref}

\makeatletter
\renewcommand{\paragraph}{%
	\@startsection{paragraph}{4}{\z@}%
	{1.5ex plus .5ex minus .2ex}{-1em}%
	{\normalfont\bfseries}}
\makeatother

\newcommand{\C}{\mathbb C}
\newcommand{\Z}{\mathbb Z}
\newcommand{\PP}{\mathbb P}
\DeclareMathOperator{\Eff}{Eff}
\newtheorem{theorem}{Theorem}[section]
\renewcommand{\thetheorem}{\ifnum\value{section}=1\relax\arabic{theorem}\else\thesection.\arabic{theorem}\fi}
\newtheorem{lemma}[theorem]{Lemma}
\theoremstyle{definition}
\newtheorem{remark}[theorem]{Remark}

\title{The moduli space of stable rational curves with seven marked points is not a Mori Dream Space}
\author{Lukas Braun}

\thanks{
This research was funded in whole or in part by the Austrian Science Fund (FWF), grants [10.55776/PAT9495823] and [10.55776/CSP1369826].
}


\begin{document}

\begin{abstract}
We construct infinitely many rigid extremal prime divisors $D_n$ on $\overline M_{0,7}$.
Our construction uses a triple cover of $\overline M_{0,7}$, ramifying over certain boundary divisors, and shows that it is birational to the Hilbert square of Vinberg's K3-surface of discriminant three. Utsumi's elliptic fibration provides us with infinitely many curves on the K3-surface, yielding divisors on the Hilbert square that finally descend to  $\overline M_{0,7}$. Our construction works in characteristic $0$ and with adaptations in every prime characteristic.
Over $\C$ we also exhibit a plane section of the pseudoeffective cone whose boundary meets an inscribed conic in a Cantor set of accumulation points of extremal divisor rays.
\end{abstract}

\maketitle

\setcounter{tocdepth}{1}
\tableofcontents
\section{Introduction}
Let $\overline M_{0,n}$ be the moduli space of stable rational curves with $n$ marked points over $\C$. In their seminal paper on Mori Dream Spaces~\cite{HuKeel2000}, Hu-Keel asked whether $\overline M_{0,n}$ is in fact a Mori dream space, which led to numerous investigations and related results since then.

Since $\overline M_{0,6}$ is a log-Fano threefold by~\cite{KeelMcKernan1996},already~\cite[Cor.~2.16]{HuKeel2000} gives a positive answer in the case $n\leq 6$. Castravet~\cite{Castravet2009} proved that indeed the sections of the $25$ boundary and $15$ \textit{Keel-Vermeire divisors} generate the Cox ring, see also~\cite{Vermeire2002,HassettTschinkel2002}. 

Non-finite generation of the Cox ring for higher $n$  was proved for $n\geq134$ by Castravet-Tevelev, for $n\geq13$ by Gonz\'alez-Karu, and for $n\geq10$ by Hausen-Keicher-Laface~\cite{CastravetTevelev2015,GonzalezKaru2016,HausenKeicherLaface2018}.

Castravet-Laface-Tevelev-Ugaglia proved the stronger property of nonpolyhedrality of the pseudoeffective cone for $n\geq10$~\cite{CastravetLafaceTevelevUgaglia2023}. Mullane obtained the same conclusion for $n\geq8$ by constructing a nonpolyhedral extremal ray of the moving cone from a stratum of meromorphic differentials~\cite[Thms. 1.2-1.4, Cor. 1.5]{Mullane2025Isoresidual}. 

There have been many further developments and attempts to get a grasp of the structure in particular of $\overline M_{0,7}$ throughout the years, of which we only can list some: Hypertree divisors and further extremal divisors provide explicit effective classes~\cite{CastravetTevelev2013,Opie2016}. Doran-Giansiracusa gave a description of the Cox ring via additive invariants, and Doran-Giansiracusa-Jensen gave a combinatorial description of its sections~\cite{DoranGiansiracusa2014,DoranGiansiracusaJensen2017}. He-Yang proved nonpolyhedrality for the blowup of $\overline M_{0,7}$ at a very general point~\cite{HeYang2021}. Computations of effective divisor classes can be found e.g. in~\cite{GonzalezGuntherZhang2020,DutourSikiricJovinelly2025}. In particular, Dutour Sikiri\'c-Jovinelly found 31 new $S_7$-orbits of extremal divisors and also exhibited Cox ring-generators whose degrees do not span extremal rays. A catalogue extending these computations led us to the cover in Section~\ref{sec:catalogue}, by which we can finally settle the question of (non-)finite generation for the Cox ring of $\overline M_{0,7}$:

\begin{theorem}\label{thm}
$\overline M_{0,7}$ is not a Mori Dream Space and it has a non-polyhedral pseudoeffective cone. In particular, there are infinitely many rigid extremal prime divisors $D_n$. This holds in characteristic $0$ and for every prime characteristic.
\end{theorem}

Consequently $\overline M_{0,n}$ is a Mori dream space exactly for $3\le n\le6$.

\paragraph{Sketch of the construction}
The sum $B$ of the $21$ boundary divisors of type $2|5$ is linearly equivalent to $-3K_X$. Taking a cube root of its defining section gives a cover $p:X'\to X$. We identify $X'$ birationally with the variety of lines on the cubic fourfold branched over the Segre cubic, and then with the Hilbert square $S^{[2]}$ of a K3 surface $S$. On $S$, the multiples of a section of infinite order give distinct rational curves $C_n$. The pairs of points that meet $C_n$ form a divisor on $S^{[2]}$. Its strict transform on $X'$ descends to the rigid extremal divisor $D_n$ on $X$. In Section~\ref{sec:cantor}, we use our construction to describe a plane section of the pseudoeffective cone over $\C$: it's boundary meets a conic in a Cantor set. 

\paragraph{AI disclosure.} All original mathematical ideas and proofs in this manuscript are the authors' own. ChatGPT 6 Astra was used only for an additional check of the arguments, proofreading and improvement of the exposition. This has not led to substantial changes, in particular not with respect to mathematical content.

The author used ChatGPT 6 Astra to produce a large amount (more than 6000 $S_7$-orbits) of generator degrees for the Cox ring in addition to those found in~\cite{DutourSikiricJovinelly2025}. Of these, 237 have $h_0>1$, with a maximum of $h_0=19$. This catalogue is not part of the present manuscript. The \textit{three first members of our series $D_n$} were among the divisors obtained.

\section{A triple cover branched along certain boundary divisors}\label{sec:catalogue}

Let $X=\overline M_{0,7}$ over $\C$. On a \textit{smooth} marked curve $(C;p_0,\ldots,p_6)$, we can choose an isomorphism $x:C\xrightarrow{\sim}\PP^1$ mapping $p_6$ to $\infty$ and (by scaling and translation) $p_0$ to $0$ and
$p_1$ to $1$. Denoting the images of $p_i$ by $x_i$, this gives the open embedding
$$
M_{0,7}\lhook\joinrel\longrightarrow\mathbb P^4,\qquad
(C;p_0,\ldots,p_6)\longmapsto[1:x_2:x_3:x_4:x_5].
$$
Thus, in homogeneous coordinates $[y_1:\cdots:y_5]$, all of $M_{0,7}$
lies in the affine chart $y_1\neq0$, where $x_i=y_i/y_1$ for
$2\le i\le5$. Its image is the open subset of $\mathbb A^4$ where
$x_2,\ldots,x_5$ are pairwise distinct and not equal to one of $0,1$.

\paragraph{Kapranov's compactification.}
Let $y_0:=0$ and write $H_{ij}:=\{y_i=y_j\}\subset\mathbb P^4$
for $0\le i<j\le5$. These fifteen hyperplanes form the complement
of $M_{0,7}$ in $\mathbb P^4$, including $H_{0j}=\{y_j=0\}$.
Kapranov's morphism $\kappa:X\to\mathbb P^4$ extends the open embedding
above. In our coordinates take the six points
$q_i:=[e_i], 1\leq i \leq 5$, $q_0:=[\sum_{i=1}^5 e_i]$, where $e_i$ are the standard basis vectors of $\C^5$.
Then $X$ is the iterated blowup at the six points, the fifteen lines spanned by pairs
and the twenty planes spanned by triples of the points, taking proper transforms in
Kapranov's order~\cite[\S\S4.2.1-4.2.4 and Theorem~4.3.3]{Kapranov1993Chow}.
Write $H$ for the pullback of the hyperplane class and $E_I$ for the
exceptional prime divisor over $\langle q_i:i\in I\rangle$, where
$I\subset\{0,\ldots,5\}$ and $1\le |I|\le3$. The classes of these form a
basis of $\operatorname{Pic}(X) \cong \Z^{42}$~\cite[\S4, immediately before Theorem~4.2]{CastravetTevelev2013}.

\paragraph{The boundary.}
For a partition $J\sqcup J^c=\{0,\ldots,6\}$ with
$|J|,|J^c|\geq2$, let $B_J=B_{J^c}$ be the divisor in $X$ of stable curves
having a node separating the marks in $J$ from those in $J^c$.
The stable curve corresponding to its general point has exactly two components. The $B_J$ are precisely the boundary
prime divisors, and their union has simple normal crossings, see e.g.~\cite[Theorem~2.7]{Knudsen1983Stacks} or~\cite[Preliminaries, Facts~1-2]{Keel1992}. In Kapranov's model, denoting by $\widetilde H_{ij}$ the strict transform of the hyperplane $ H_{ij}$, they are
$$
\begin{aligned}
	B_{I\cup\{6\}}&=E_I, &&1\le |I|\le3;\\
	B_{\{i,j\}}&=\widetilde H_{ij}, &&0\le i<j\le5.
\end{aligned}
$$
We write $B_{ij}$ for $B_{\{i,j\}}$ and define
$$
B:=\sum_{0\le i<j\le6}B_{ij}
=\sum_{0\le i<j\le5}\widetilde H_{ij}+\sum_{i=0}^5E_{\{i\}} .
$$
This is the reduced sum of the 21 boundary divisors of type $2|5$.
Note that there are another 35 boundary divisors of type $3|4$.

\paragraph{The section cutting out $B$.}
The product of the fifteen equations corresponding to the collision of points is the homogeneous polynomial
$$
\Delta(y)=\prod_{1\le i<j\le5}(y_j-y_i) \cdot \prod_{1\le k\le5}y_k.
$$
Any blowup center $\langle q_i:i\in I\rangle$ lies in precisely
$\binom{6-|I|}{2}$ collision hyperplanes, one for each pair of labels
outside $I$. Thus the total transform of $A:=\{\Delta=0\}$ has
multiplicities $10,6,3$ along the exceptional divisors over points,
lines and planes, respectively, and is linearly equivalent to $15H$. Writing
$E^{(k)}=\sum_{|I|=k}E_I$, this gives
\begin{equation}
	B \sim 15H-9E^{(1)}-6E^{(2)}-3E^{(3)} \sim -3K_X .                 \label{eq:K}
\end{equation}
Note that the last equality comes from the canonical divisor formula for a blowup in a
smooth centre of codimension $c$. Castravet-Tevelev use $B \sim -3K_X$ in their proof that the
Mori cone for $n=7$ is generated by $F$-curves~\cite[proof of Corollary~2.9]{CastravetTevelev2012Rigid}. Moon uses $B$ together with the reduced sum of the boundary divisors of type $3|4$ to investigate the projective models of $X$ defined by symmetric effective divisors~\cite[Corollary~2.6 and Theorem~1.1]{Moon2017Seven}.

To identify the section cutting out B itself, set $K=\mathbb C(X)=\mathbb C(x_2,x_3,x_4,x_5)$.
On our chart $\{y_1\neq0\}$, with $x_0=0,\ x_1=1$, define
\begin{equation}
	\delta:=\frac{\Delta(y)}{y_1^{15}}
	=\prod_{i=2}^5 x_i(x_i-1)\prod_{2\le i<j\le5}(x_j-x_i),
	\qquad \theta:=dx_2\wedge dx_3\wedge dx_4\wedge dx_5 .        \label{eq:pair-product}
\end{equation}
On $\mathbb P^4$, $\theta$ has divisor $-5H_{01}$;
pulling back gives
$$
K_\theta:=\operatorname{div}_X(\theta)
=-5\kappa^*H_{01}+3E^{(1)}+2E^{(2)}+E^{(3)} .
$$
The multiplicities computed above now give the equality of \emph{divisors}
$$
\operatorname{div}_X(\delta)=B+3K_\theta .
$$
Consequently, with $\mathcal A=\mathcal O_X(-K_X)$, the section
$$
s_B:=\delta\,(\theta^{-1})^{\otimes3}
$$
of $\mathcal A^{\otimes3}$ is regular and has zero divisor exactly $B$.
In particular, $\delta$ is the coefficient of $s_B$ in the frame
$(\theta^{-1})^{\otimes3}$ on $M_{0,7}$, where it is nowhere zero.

\paragraph{The cover.}
We have just seen that the effective Cartier divisor $B$ has a cube root $\mathcal A$ of its line bundle.
We denote the fiber over $x\in X$ of the total space of $\mathcal A$ by $\mathcal A_x$. In this total space set
$$
X'=\{(x,v):v\in\mathcal A_x,\ v^{\otimes3}=s_B(x)\},
\qquad p:X'\longrightarrow X,\quad (x,v)\longmapsto x .
$$

\begin{lemma}\label{lem:pair-cubic-cover}
	The morphism $p:X'\to X$ is a finite flat cyclic cover of degree three.
	Over $M_{0,7}$, the frame $\theta^{-1}$ identifies the total space of $\mathcal A$ with $M_{0,7}\times\mathbb A^1_\tau$ via $v=\tau\theta^{-1}$. In this trivialization $X'$ is given by $\tau^3=\delta$, and $p$ is just the projection to the first factor.
	On $X=\overline M_{0,7}$, the branch locus is exactly
	$B=\sum_{0\le i<j\le6}B_{ij}$. There is a unique prime divisor
	$R_{ij}$ above each $B_{ij}$, and $p^*B_{ij}=3R_{ij}$.
	Moreover, $X'$ is the normalization of $X$ in
	$L=\mathbb C(x_2,x_3,x_4,x_5)(\tau)$.
\end{lemma}

\begin{proof}
	By our above computations, the properties of the construction follow immediately.
	In particular, for every pair $i,j$, the section $s_B$ vanishes to order one
	along $B_{ij}$. The cyclic-cover ramification formula therefore
	gives index three along each of these divisors, and the cover is
	\'etale over $X\setminus B$~\cite[Lemma~3.15]{EsnaultViehweg1992}.
	Since the cover has degree three, there is just one prime $R_{ij}$
	above $B_{ij}$, with $p^*B_{ij}=3R_{ij}$.
	Finally, since $B$ is reduced, the cyclic-cover construction
	identifies $X'$ with the normalization of $X$ in $L$ and gives
	$p_*\mathcal O_{X'}=\mathcal O_X\oplus\mathcal A^{-1}
	\oplus\mathcal A^{-2}$, hence finite flatness of degree three~\cite[\S3.7]{EsnaultViehweg1992}.
\end{proof}

We write $R=\sum_{i<j}R_{ij}$ as well as $\eta:=s_B^{-1}=\theta^{\otimes3}/\delta$ on $X$ and $\Theta=p^*\theta/\tau$ on $X'$. By this, since  a local parameter pulls back with order three along $R_{ij}$,
and its differential contributes order two in the pullback of a top
form, we get:
\begin{equation}
	\Theta^{\otimes3}=p^*\eta,\qquad
	\operatorname{div}_{X'}(\Theta)=R .                          \label{eq:tensor}
\end{equation}

\begin{remark}
In~\cite[\S 4]{Moon2017Seven}, Moon constructs the model $Y=\overline M_{0,7}^{\,1}$ by two small modifications
of $X$, followed by contracting the strict transforms
of the divisors $B_{ij}$.
Let $Y'$ be the normalization of $Y$ in $L$, which fits in the diagram
$$
\begin{tikzcd}[column sep=large,row sep=large]
 X' \arrow[r,dashed] \arrow[d,"p"'] & Y' \arrow[d,"q"] \\
 X \arrow[r,dashed] & Y
\end{tikzcd}
$$
Take $K_Y=\operatorname{div}_Y(\theta)$. On $Y$ we have
$\operatorname{div}_Y(\delta)=3K_Y$. If $K_Y=\operatorname{div}_Y(h)$, the function $\delta/h^3$ would have trivial divisor and would be constant on the normal projective variety $Y$. Since every nonzero complex constant is a cube, this would make $\delta$ a cube in $\C(X)$, contrary to $[L:K]=3$. Thus $K_Y$ has order three in the divisor class group. Thus $q:Y'\to Y$ is the global
canonical cover associated with $3K_Y\sim0$.
\end{remark}

\section{The triple cover is birational to a variety of lines}\label{sec:field}

Forgetting the seventh point $p_6$ gives a
morphism $M_{0,7}\to M_{0,6}$. Considering the remaining six points
up to all projective coordinate changes of the line, we get the equal-weight
GIT compactification of $M_{0,6}$: the Segre cubic threefold
$\Sigma=(\mathbb P^1)^6/\!/\mathrm{SL}_2$. For its classical
description by six-point invariants see~\cite[\S35]{Coble1929Theta}.
The composition map $M_{0,7}\to M_{0,6}\hookrightarrow\Sigma$ extends to
$X\to\overline M_{0,6}\to\Sigma$, with stabilization after forgetting
$p_6$ on the boundary~\cite[\S2, diagram~(2.2), and \S3]{Dolgachev2016Segre}.

We describe this map to the Segre cubic in the coordinates
of Section~\ref{sec:catalogue}. Write the six marked points
$p_0,\ldots,p_5$ on their line as $p_i=[\alpha_i:\beta_i]$ and define
the \emph{bracket}
$$
[ij]:=\alpha_i\beta_j-\alpha_j\beta_i .
$$
In the coordinate $x$ sending $p_0,p_1,p_6$ to $0,1,\infty$, taking
$(\alpha_i,\beta_i)=(x_i,1)$ gives $[ij]=x_i-x_j$, with
$x_0=0$ and $x_1=1$.

A matching $m$ is a partition of $\{0,\ldots,5\}$ into three pairs.
Its bracket product is
$$
P_m:=\prod_{ij\in m}[ij],
$$
where we take $i<j$ in each pair, see e.g.~\cite[\S35]{Coble1929Theta} or~\cite[\S2]{HowardMillsonSnowdenVakil2009}.
Rescaling the representative of $p_i$ multiplies every $P_m$ by
the same scalar, since each matching uses $i$ exactly once.
Likewise, a common coordinate change given by $g\in\mathrm{GL}_2$
multiplies every $P_m$ by $(\det g)^3$. Thus the projective ratios
of the matching products depend only on the six-pointed curve.
We use the five matchings
$$
a=01|23|45,\quad b=01|24|35,\quad c=02|13|45,\quad
d=02|14|35,\quad h=03|14|25.
$$
Their bracket products $P_a,\ldots,P_h$ form a basis of the span of all fifteen matching
products $P_m$. The others can be expressed in this basis using the
Pl\"ucker relation $[ij][kl]-[ik][jl]+[il][jk]=0$. For example,
$$
P_{05|14|23}=[14]\bigl([03][25]-[02][35]\bigr)=P_h-P_d.
$$

In the homogeneous coordinates $[y_1:\cdots:y_5]$ on the source, with $y_0=0$, the forgetful map $M_{0,7} \to \Sigma$ is
$$
[y_1:\cdots:y_5]\longmapsto
[P_a(y):P_b(y):P_c(y):P_d(y):P_h(y)]
\in\PP^4_{[a:b:c:d:h]}.
$$
On our chart $y_1\neq0$ we can just write $P=P(x)=(P_a(x),P_b(x),P_c(x),P_d(x),P_h(x))$.
In these target coordinates the Segre cubic has equation
\begin{equation}\label{eq:segre-match}
	F(a,b,c,d,h):=h(bc-ad)-d(a-b)(a-c)=0.
\end{equation}
Indeed, after forgetting $p_6$, choose the coordinate $u$ sending
$p_0,p_1,p_5$ to $0,1,\infty$, and put $u_i=u(p_i)$ for $i=2,3,4$.
Using $(u_i,1)$ for the finite points and $(1,0)$ for $p_5$, the
same bracket products give
$$
[a:b:c:d:h]
=[u_2-u_3:u_2-u_4:u_2(1-u_3):u_2(1-u_4):u_3(1-u_4)],
$$
on which $F$ vanishes. Conversely, on a dense open set we recover
the three parameters by
\begin{equation}\label{eq:segre-inverse}
	u_2=\frac{c-d}{a-b},\qquad
	u_3=\frac{u_2h}{d},\qquad
	u_4=1-\frac{d(1-u_3)}{c}.
\end{equation}
Thus the image has dimension three. Since $F$ is linear in $h$
with coprime coefficients $bc-ad$ and $-d(a-b)(a-c)$, it is
irreducible, and its zero locus is the closure of the image.

The Segre cubic $\Sigma$ has ten nodes~\cite[Proposition~2.2]{Dolgachev2016Segre}, which are indexed by unordered partitions
$I\sqcup I^c=\{0,\ldots,5\}$ into two triples, corresponding to the degenerate point  configuration
with three marks at $0$ and three at $1$. To obtain the node
corresponding to $I\mid I^c$, put $\epsilon_i=0$ for $i\in I$ and
$\epsilon_i=1$ for $i\in I^c$. Its coordinates are then
$$
[a:b:c:d:h]=[P_a(\epsilon):P_b(\epsilon):P_c(\epsilon):
P_d(\epsilon):P_h(\epsilon)].
$$

\paragraph{The cubic fourfold and its lines.}
In $\mathbb P^5$ with coordinates $[a:b:c:d:h:z]$, let
$$
W=\{z^3=F(a,b,c,d,h)\}.
$$
Projection onto the first five coordinates makes $W$ the cyclic
triple cover of $\mathbb P^4$ branched along $\Sigma$~\cite[\S2, p.~3]{Mukai2022Segre}.
A point of $W$ is singular exactly when $z=0$ and the five partial
derivatives of $F$ vanish. Thus $W$ has one singular point above
each of the ten nodes of $\Sigma$. Each has multiplicity two,
since $z^3-F=0$ has the same nonzero quadratic
part as $-F$.

A point of $\operatorname{Gr}(2,6)$ is a two-dimensional subspace
$U\subset\mathbb C^6$ and determines the line
$\ell=\mathbb P(U)\subset\mathbb P^5$. Write $F_1(W)$ for the
locus of those $U$ for which $\ell\subset W$. On the Grassmannian
chart where $a$ and $z$ restrict to independent linear forms on $U$,
we may use $s=a$ and $t=z$ as coordinates on $U$. The other four
coordinates are then linear forms in $s,t$. Thus the line has a
unique parametrization of the form
$$
[a:b:c:d:h:z]=[s:su+tv:t],\qquad [s:t]\in\mathbb P^1,
$$
where the middle block $su+tv$ gives the four coordinates
$(b,c,d,h)$. The entries of $u,v\in\mathbb C^4$ are the eight
affine coordinates on this Grassmannian chart.

In the chart, the line $\ell=[s:su+tv:t]$ lies on $W$
exactly when $F(s,su+tv)=t^3$ for all $[s:t]\in\PP^1$. Both sides
are binary cubics, so it suffices to test equality at
$[1:0],[0:1],[1:1],[1:-1]$. Thus the equations for $F_1(W)$ on
this chart are
$$
F(1,u)=0,\quad F(0,v)=1,\quad
F(1,u+v)=1,\quad F(1,u-v)=-1.
$$
Globally, $F_1(W)$ is irreducible of dimension four by~\cite[Proposition~3.5 and its proof]{Lehn2015Starr} since the singularities
of $W$ are isolated and have multiplicity two.

Projection to $\PP^4$ maps every line $\ell\subset W$
isomorphically onto a line $\lambda\subset\PP^4$, since the center lies outside $W$. 
If $z$ vanishes identically on $\ell$, then
$\lambda\subset\Sigma$ and $\ell$ is its unique lift. Otherwise
$z$ has one simple zero on $\ell$. Since $F=z^3$ along $\ell$,
the projected line $\lambda$ meets $\Sigma$ at the image of that
point with multiplicity three.

Conversely, suppose $\lambda$ meets $\Sigma$ only at $[P]$, with
multiplicity three. Write $\lambda=[sP+tV]$, so that $t=0$ at $[P]$.
Then $F(sP+tV)=c\,t^3$ for some $c\in\C^*$, and the three lines
in $W$ above $\lambda$ are
$$
[s:t]\longmapsto[sP+tV:\tau t],\qquad \tau^3=c.
$$
The next theorem constructs $\lambda$ from a smooth seven-pointed
curve. The first six marks give the contact point $[P]$, and the
seventh determines the direction of $\lambda$. In the coordinates
of Section~\ref{sec:catalogue}, we obtain $c=\delta$. Thus choosing
a lift of $\lambda$ amounts to choosing a point of $X'$ above the
marked curve.

\begin{theorem}\label{thm:lines-cover}
	For a matching $m$, define
	$$
	\sigma_m:=\sum_{ij\in m}x_i x_j,\qquad V_m:=P_m(x)\sigma_m,\qquad V:=(V_m)_m.
	$$
	For every $(x_2,x_3,x_4,x_5,\tau)\in p^{-1}(M_{0,7})$,
	with $\tau^3=\delta(x)$, the line
	\begin{equation}\label{eq:matching-line}
		\ell(x,\tau)=
		\bigl\{[sP(x)+tV:t\tau]\mid[s:t]\in\PP^1\bigr\}
	\end{equation}
	is contained in $W$. The assignment $(x,\tau)\mapsto\ell(x,\tau)$
	defines a morphism $p^{-1}(M_{0,7})\to F_1(W)$ and induces a
	birational equivalence $X'\dashrightarrow F_1(W)$.
\end{theorem}

\begin{proof}
	We first prove the identity
	\begin{equation}\label{eq:flex-identity}
		F(sP+tV)=t^3\delta.
	\end{equation}
	Use the homogeneous coordinates of Section~\ref{sec:catalogue},
	with $y_0=0$, and evaluate $P_m$ and $V_m$ at
	$y_0,\ldots,y_5$. We will show that every coefficient of
	$F(sP+tV)$ is divisible by
	$$
	\Delta(y)=\prod_{0\le i<j\le5}(y_j-y_i).
	$$
	
	For any matching $m$, write $L_m(a,b,c,d,h)$ for the linear form
	expressing $P_m$ in our chosen basis
	$P_a,P_b,P_c,P_d,P_h$. Thus $L_m(P)=P_m$. We claim that also
	$L_m(V)=P_m\sigma_m$ for every matching $m$. To see this, evaluate the bracket products
	on the six points $y_k^r:=[y_k:1-r y_k^2]\in\PP^1$. The bracket for the
	pair $k,l$ is
	$$
	y_k(1-r y_l^2)-y_l(1-r y_k^2)
	=(y_k-y_l)(1+r y_k y_l).
	$$
	Hence the product for $m$ is
	$$
	P_m(y_0^r,\ldots,y_5^r)=P_m(y)(1+r\sigma_m)+O(r^2). 
	$$
	For the five basis matchings this gives
	the vector $P+rV+O(r^2)$. Applying $L_m$ and comparing
	coefficients of $r$ proves our claim.
	
	Now fix $i,j$ and write $k,l,u,v$ for the remaining four labels.
	The closure of the image of the hyperplane $y_i=y_j$ in the source $\PP^4$ - corresponding to the collision of the marked points $p_i$ and $p_j$ - under our  map is the plane
	$$
	\Pi_{ij}
	=\{L_{ij|kl|uv}=L_{ij|ku|lv}=0\}
	\subset \PP^4_{[a:b:c:d:h]},
	$$
	which lies on the Segre cubic $\Sigma$.
	Indeed, both bracket products contain the factor $y_i-y_j$, and
	the third matching containing $ij$ is a linear combination of these
	two by the Pl\"ucker relation.
For $ij=01$ the two equations may be taken to be $a=b=0$, and $F(0,0,c,d,h)=0$. Permuting the six labels acts linearly on the matching products and preserves their image $\Sigma$. Thus every plane $\Pi_{ij}$ is contained in $\Sigma$. 
	When $y_i=y_j$, the two forms displayed above vanish on $P$ and on $V$. They therefore vanish on $sP+tV$,
	so $F(sP+tV)=0$ whenever $y_i=y_j$.
	
	Thus all fifteen distinct factors $y_j-y_i$ divide each coefficient
	of $F(sP+tV)$. Since $P_m$ and $V_m$ have degrees $3$ and $5$,
	the coefficients of $s^3,s^2t,st^2,t^3$ have degrees
	$9,11,13,15$, respectively. The first three therefore vanish, and
	the last is a constant multiple of $\Delta$. Evaluating at a fixed $y$ in fact  gives $F(V)=\Delta$. Setting $y_1=1$ now gives
	\eqref{eq:flex-identity} in our original coordinates
	$x_0=0,\ x_1=1$.
	Since $(t\tau)^3=t^3\delta$, the line \eqref{eq:matching-line}
	lies in $W$.
	
	Now, vice versa, we recover the marked curve and $\tau$ from the line $\ell$.
	Along $\ell$ we have $a=P_a s+V_a t$ and $z=\tau t$.
	Since the six finite marks are distinct, both $P_a$ and $\tau$
	are nonzero. Thus the change of coordinates from $[s:t]$ to
	$[a:z]$ has nonzero determinant $P_a\tau$, and $[a:z]$ are
	homogeneous coordinates on $\ell$. Its intersections with
	$z=0$ and $a=0$ therefore have unique representatives of the form
	$$
	A:=(1,A_b,A_c,A_d,A_h,0),\qquad
	B:=(0,B_b,B_c,B_d,B_h,1),
	$$
	respectively. By \eqref{eq:matching-line} and $V_j=P_j\sigma_j$,
	$$
	A_j=\frac{P_j}{P_a},\qquad
	B_j=\frac{P_j(\sigma_j-\sigma_a)}{\tau}
	\qquad \mathrm{for~}j=b,c,d,h.
	$$
	The point $A$ thus gives $[P]$. Hence \eqref{eq:segre-inverse}
	recovers the coordinates $u_2,u_3,u_4$ of the six-pointed curve
	in the coordinate $u$ sending $p_0,p_1,p_5$ to $0,1,\infty$.
	
	To recover the seventh mark, put $r_j=B_j/A_j$. Then
	$r_j=P_a(\sigma_j-\sigma_a)/\tau$, so the line determines
	$$
	k:=\frac{r_b}{r_d-r_c}
	=\frac{\sigma_b-\sigma_a}{\sigma_d-\sigma_c}
	=\frac{(x_2-x_5)(x_4-x_3)}{(1-x_5)(x_4-x_3)}
	=\frac{x_2-x_5}{1-x_5}.
	$$
	In the original coordinate $x$, where $p_0,p_1,p_6$ have
	coordinates $0,1,\infty$, we have
	$u=(1-x_5)x/(x-x_5)$. In particular $u_2=x_2/k$, so
	$$
	x_2=ku_2,\qquad x_5=\frac{x_2-k}{1-k},\qquad
	x_i=\frac{u_i x_5}{u_i-1+x_5}\quad(i=3,4).
	$$
	This recovers all four coordinates on $M_{0,7}$. Finally,
	evaluating $P_a,\sigma_a,\sigma_b$ at the above coordinates
	gives
	$$
	\tau=\frac{P_a(\sigma_b-\sigma_a)}{r_b}.
	$$
	
	These formulas give a rational inverse on the image. The image
	therefore has dimension four and is dense in the irreducible
	fourfold $F_1(W)$. This proves the birational equivalence.
\end{proof}

\section{The variety of lines is birational to a Hilbert square}\label{sec:hilbert}

Mukai~\cite{Mukai2022Segre} identifies the variety of lines
$F_1(W)$ birationally with the Hilbert square
of Vinberg's K3 surface $X_3$ of discriminant three~\cite{Vinberg1983}.
We describe the surface and the birational map in the coordinates of
Section~\ref{sec:field}.

\paragraph{The K3 surface of lines through a singular point.}
On $W$, we choose the singular point $o=[0:0:0:0:1:0]$ lying above the node
$012\mid345$ of the Segre cubic. 
Any line through $o$ has points
$[t a:t b:t c:t d:h:t z]$, with $[t:h]\in\mathbb P^1$, so it is given by the direction $[a:b:c:d:z] \in \mathbb P^4$. Put
$$
q=bc-ad,\qquad g=z^3+d(a-b)(a-c).
$$
The equation $z^3=F(a,b,c,d,h)$ of $W$ equals $h q=g$.
Substituting the parametrization of the line gives $t^2(hq-tg)=0$. Thus the entire line lies
on $W$ exactly when $q=g=0$ holds for the direction $[a:b:c:d:z]$. The surface parametrizing the lines through $o$ lying in $W$ thus is
$$
S_0=\{bc-ad=z^3+d(a-b)(a-c)=0\}\subset\mathbb P^4_{[a:b:c:d:z]} .
$$

Forgetting $z$ gives a triple cover of the smooth quadric
$Q=\{bc=ad\}\subset\mathbb P^3_{[a:b:c:d]}$. On the chart $d\neq0$ put
$$
x=b/d,\qquad y=c/d,\qquad \xi=z/d.
$$
Since $a/d=xy$, a point of $S_0$ in this chart has coordinates
$[xy:x:y:1:\xi]$. Substitution into $z^3=-d(a-b)(a-c)$ gives
\begin{equation}\label{eq:quadric-cover}
	\xi^3=-x(x-1)y(y-1).
\end{equation}
Here $x,y$ are affine coordinates on the two factors of
$Q\simeq\mathbb P^1\times\mathbb P^1$. The branch divisor is the union of the fibres $x=0,1,\infty$ and $y=0,1,\infty$. They meet at the nine points $(x,y)\in\{0,1,\infty\}^2$
of $Q$. Above each of these points lies a unique point of $S_0$, with
local equation $w^3=uv$. These nine $A_2$ singularities are the only
singular points of $S_0$.

We can now identify the surface in Utsumi's coordinates, see~\cite[\S2]{Utsumi2016}. With $\gamma^3=-16$, let
$$
r_1=2x-1,\qquad r_2=2y-1,\qquad Z=\gamma\xi.
$$
Then \eqref{eq:quadric-cover} reads
$$
Z^3=(r_1^2-1)(r_2^2-1),
$$
which is precisely Utsumi's affine model of $X_3$.
Thus the minimal resolution $\rho:S\to S_0$ is birational to $X_3$.
In fact $S\simeq X_3$: adjunction for the complete intersection
$S_0\subset\mathbb P^4$ gives $\omega_{S_0}\simeq\mathcal O_{S_0}$,
and the resolution of its $A_2$ singularities is crepant. Hence $S$
is already the smooth minimal model. Over each of the nine singular
points, $\rho$ has two exceptional $(-2)$-curves.

\paragraph{The residual line.}
The Hilbert square $S^{[2]}=\operatorname{Hilb}^2(S)$ parametrizes
length-two subschemes of $S$. Its dense open subset of reduced
subschemes parametrizes unordered pairs of distinct points of $S$.
A length-two subscheme supported at a single point $p\in S$ is
uniquely determined by the choice of a tangent line $L\subset T_pS$. Since $S$ is a projective K3 surface,
$S^{[2]}$ is a smooth projective holomorphic symplectic fourfold, see~\cite[Theorem~2.4]{Fogarty1968} and~\cite[\S6, Theorem~3]{Beauville1983}.
In particular, its canonical bundle is trivial.

\begin{theorem}\label{thm:hilbert-residual}
	For a general pair of points of $S$, the two corresponding lines
	$L_1,L_2\subset W$ through $o$ span a plane $\Pi$. The plane cubic
	$\Pi\cap W$ is the union of $L_1$, $L_2$, and a third line $\ell$,
	which we call the \textit{residual line}. Sending the pair to $\ell$ defines a birational	map $\psi:S^{[2]}\dashrightarrow F_1(W)$, equivariant for the action
	induced by $z\mapsto\zeta z$, where $\zeta$ is a primitive cube
	root of unity.
\end{theorem}

\begin{proof}
	We describe $\psi$ and its inverse in coordinates, following~\cite[proof of Theorem~3.6]{Lehn2015Starr}. For two general points of $S$, let $m\subset\mathbb P^4_{[a:b:c:d:z]}$
	be the line joining their images on $S_0$. In homogeneous coordinates
	$[\lambda:\mu]$ on $m$ with the two points at $[1:0]$ and
	$[0:1]$, both $q|_m$ and $g|_m$ vanish at the two points. Hence
	$$
	q|_m=A\lambda\mu,\qquad
	g|_m=\lambda\mu(B_1\lambda+B_2\mu).
	$$
	For general points $A\neq0$, since the line joining their images on $Q$
	meets $Q$ only at those two points. The plane over $m$ through $o$
	has coordinates $[\lambda:\mu:h]$, and its section by $W$ is
	$$
	\lambda\mu\bigl(Ah-B_1\lambda-B_2\mu\bigr)=0 .
	$$
	The factors $\lambda,\mu$ give the two lines through $o$.
	The residual line $\ell$ has equation $Ah=B_1\lambda+B_2\mu$.
	
	Conversely, a line $\ell\subset W$ missing $o$ projects isomorphically
	to a line $m$ in the space of directions. On $m$ the coordinate $h$
	of its lift is a linear form, so $g|_m=h\,q|_m$. If $q|_m\neq0$, the quadratic $q|_m$ cuts out the length-two
	subscheme $m\cap S_0$, since $S_0=\{q=g=0\}$. For the line $\ell$
	constructed above, this intersection consists of the original two points. This
	recovers a general pair on $S$, so the image has dimension four and is
	dense in the irreducible fourfold $F_1(W)$ of Section~\ref{sec:field}.
	The automorphism $z\mapsto\zeta z$ fixes $o$, preserves $q,g$,
	and lifts to the minimal resolution. It therefore commutes with both
	constructions.
\end{proof}

In particular, for the line $\ell(x,\tau)=[sP+tV:t\tau]$ of
Section~\ref{sec:field}, the pair on $S$ is obtained by solving the
single binary quadratic
\begin{equation}\label{eq:pair-quadratic}
	(sP_b+tV_b)(sP_c+tV_c)-(sP_a+tV_a)(sP_d+tV_d)=0 .
\end{equation}
At each root one takes
$[sP_a+tV_a:sP_b+tV_b:sP_c+tV_c:sP_d+tV_d:t\tau]$ on $S_0$
and lifts it by $\rho$. This description applies on the open set where
the roots are distinct and the two points are smooth. Here, as in
Section~\ref{sec:field}, $P,V$ are evaluated at
$(0,1,x_2,x_3,x_4,x_5)$ and $\tau^3=\delta(x)$.

\paragraph{The birational maps and their exceptional divisors.}
Let $\alpha:X'\dashrightarrow F_1(W)$ be the map of
Section~\ref{sec:field} and put $\beta=\alpha^{-1}\circ\psi$.
The action of $G=\langle\zeta\rangle$ on $X'$ is
$\tau\mapsto\zeta\tau$. The formula for $\ell(x,\tau)$ shows that
$\alpha$ is equivariant. Taking quotients gives the following commutative
diagram:
$$
\begin{tikzcd}[column sep=large,row sep=large]
	Y' \arrow[d] &
	X' \arrow[l,dashed] \arrow[r,dashed,"\alpha"] \arrow[d,"p"] &
	F_1(W) \arrow[d] &
	S^{[2]} \arrow[l,dashed,"\psi"'] \arrow[d] \\
	Y & X \arrow[l,dashed] \arrow[r,dashed] &
	F_1(W)/G & S^{[2]}/G \arrow[l,dashed]
\end{tikzcd}
$$
Here $Y$ is Moon's model and $Y'$ its normalization in
$\mathbb C(X')$, as in Section~\ref{sec:catalogue}. All horizontal
arrows are birational and all vertical arrows are finite quotients.

For an irreducible curve $C\subset S$, write $I_C\subset S^{[2]}$
for the closure of the pairs of distinct points with one point on
$C$. It is a prime divisor, as it is the closure of the image of a
dense open subset of $C\times S$ with generically finite fibres.
Write $f_1,f_2:S\to\mathbb P^1$ for the two projections through $Q$,
and $T_j$ for the closure of the pairs of distinct points on a common
smooth fibre of $f_j$.

\begin{remark}\label{rem:divisors-models}
	The divisors on $S^{[2]}$ contracted by $\psi$
	are the eighteen $I_E$ for the exceptional curves of $S\to S_0$
	and the two divisors $T_1,T_2$ from the rulings.
	Every prime divisor on $S^{[2]}$, however, has a prime strict transform
	under $\beta:S^{[2]}\dashrightarrow X'$. The only primes contracted
	in the reverse direction are the $21$ ramification divisors $R_{ij}$,
	as we show in Lemma~\ref{lem:section-divisors}.
	In particular, $Y'$ and $S^{[2]}$ are isomorphic in codimension one.
	Indeed, Moon's map $X\dashrightarrow Y$ contracts precisely the pair
	boundaries and its inverse contracts no divisor~\cite[Theorem~1.1 and \S4.6]{Moon2017Seven}.
	Normalizing in the common cubic extension therefore removes precisely
	the $R_{ij}$ from $X'$ in codimension one.
\end{remark}

\section{Elliptic sections and incidence divisors on the triple cover}
\label{sec:elliptic-incidence}

On the affine chart of $S_0$ from Section~\ref{sec:hilbert}, the
surface equation and the generator $g$ of the deck group $G$ are
$$
Z^3=(r_1^2-1)(r_2^2-1),\qquad
g(r_1,r_2,Z)=(r_1,r_2,\zeta Z).
$$
We use the same symbols for the pullbacks of $r_1,r_2,Z$ to the
minimal resolution $\rho:S\to S_0$ and for the lift of $g$ to $S$.
We seek infinitely many curves on $S$ preserved by $g$.
The two projections to $r_1$ and $r_2$ give elliptic fibrations
on $S$, each with only three sections~\cite[Table~1]{Utsumi2016}.
We therefore use another fibration by Utsumi, whose sections include a point
of infinite order on the generic fibre.

\paragraph{The elliptic fibration on $S$.}
Utsumi's 'Fibration~3`, see~\cite[\S 4]{Utsumi2016}, is the morphism $f:S\to\mathbb P^1_u$ for which
the pullback of the affine coordinate $u$ on the base is
$$
f^*u=\frac{Z}{r_1^2-1}
\qquad\text{in }\mathbb C(S).
$$
Writing $u$ also for $f^*u$, set $k=u^3$.
Substitution of $Z=u(r_1^2-1)$ into the surface equation gives
the affine equation of the generic fiber of $f$:
\begin{equation}\label{eq:elliptic-quartic}
	r_2^2=1+k(r_1^2-1)^2 .
\end{equation}
Utsumi's change of variables identifies the smooth projective model of
the quartic \eqref{eq:elliptic-quartic} with the elliptic curve
$\mathcal E$ given in affine Weierstrass coordinates by
\begin{equation}\label{eq:weier}
	Y_e^2=X_e^3+4kX_e^2-4kX_e .
\end{equation}
Under this isomorphism, the point $(r_1,r_2)=(1,1)$ on the quartic
maps to the point at infinity $[0:1:0]$ on the projective
Weierstrass cubic. We use this point as the identity for the group
law on $\mathcal E$.

To express a point $(X_e,Y_e)$ of $\mathcal E$ in the surface
coordinates used above, we will use the inverse of Utsumi's change
of variables. On the open set $Y_e\neq4k$ and with $T=\frac{2(X_e+4k)}{Y_e-4k}$, it is given by
$$
r_1=1+T,\qquad r_2=\frac{X_eT^2}{2}-1,\qquad
Z=u(r_1^2-1).
$$

\paragraph{The section curves.}
Let $\operatorname{MW}(f)$ denote the group of sections of
$f:S\to\mathbb P^1_u$. Restriction to the generic fibre identifies
$\operatorname{MW}(f)$ with the group of $\mathbb C(u)$-rational
points of the elliptic curve $\mathcal E$ in \eqref{eq:weier}.
Utsumi shows that the point $P=(X_e,Y_e)=(1,-1)$ has infinite order.

For $n\geq1$, let $s_n:\mathbb P^1_u\to S$ be the section
corresponding to $nP$, and let $C_n=s_n(\mathbb P^1_u)\subset S$.
The curves $C_n$ are pairwise distinct because the points $nP$
are distinct on the generic fibre. The restriction of $f$ to $C_n$
is an isomorphism onto $\mathbb P^1_u$, with inverse $s_n$.
Thus $C_n$ is a smooth rational curve and $C_n^2=-2$ by adjunction
on the K3 surface $S$.

The deck transformation $g$ sends $u=Z/(r_1^2-1)$ to $\zeta u$
and fixes $k=u^3$. It also fixes $X_e,Y_e$, since Utsumi's
formulas express them in terms of $r_1,r_2,k$, all fixed by $g$~\cite[\S4, (4.5)]{Utsumi2016}. Both the Weierstrass equation and
$P=(1,-1)$ are defined over $\mathbb C(k)$. Hence the coordinates
of every multiple $nP$ are rational functions of $k$, so they
take the same values at $u$ and $\zeta u$. It follows that
$g(s_n(u))=s_n(\zeta u)$, and therefore $g(C_n)=C_n$.

\paragraph{Incidence on the Hilbert square.}
For each section curve $C_n\subset S$, let $I_n\subset S^{[2]}$
be the reduced locus of length-two
subschemes whose support meets $C_n$. For a pair of distinct
points, this means that at least one of the two lies on $C_n$.
The pairs with exactly one point on $C_n$ form a dense open
subset of $I_n$, isomorphic to $C_n\times(S\setminus C_n)$
via $(a,q)\mapsto\{a,q\}$. Hence $I_n$ is a prime divisor. In particular, it is the divisor denoted $Z_C$ for a smooth rational curve $C\subset S$ in~\cite[Theorem~4.2]{AmerikVerbitsky2022}. Since $g(C_n)=C_n$, it is also $G$-invariant.

Fix $q\in S\setminus C_n$. The morphism
$$
\iota_q:S\setminus\{q\}\longrightarrow S^{[2]},\qquad
a\longmapsto\{a,q\},
$$
embeds $C_n$ as a curve $\Gamma_{n,q}:=\iota_q(C_n)\subset I_n$.
As $q$ varies in $S\setminus C_n$, these curves cover the dense
open subset of $I_n$ described above.

To compute $I_j\cdot\Gamma_{n,q}$, choose $q$ outside
$C_n\cup C_j$. Then, since the fixed point $q$ lies
outside $C_j$, a local equation for the incidence condition on
$\{a,q\}$ is simply a local equation for $a\in C_j$. So $\iota_q^*I_j=C_j$ as divisors on $S\setminus\{q\}$. Hence
$$
I_j\cdot\Gamma_{n,q}
=\deg\!\left(\mathcal O_S(C_j)|_{C_n}\right)
=C_j\cdot C_n .
$$
For $j=n$ this gives $I_n\cdot\Gamma_{n,q}=C_n^2=-2$.
Since the curves $\Gamma_{n,q}$ cover a dense open subset of $I_n$,
the moving-curve criterion~\cite[Lemma~4.1]{ChenCoskun2014} shows that
$I_n$ spans an extremal ray of the pseudoeffective cone of $S^{[2]}$
and that
\begin{equation}\label{eq:rigid}
	|mI_n|=\{mI_n\}\qquad(m\geq1).
\end{equation}
For $j\neq n$, choose $q\notin C_n\cup C_j$. If $[I_j]$ and
$[I_n]$ spanned the same ray, we would have $[I_j]=\lambda[I_n]$
for some $\lambda>0$. Intersecting with $\Gamma_{n,q}$ would give
$$
0\le C_j\cdot C_n
= I_j\cdot\Gamma_{n,q}
= \lambda\, I_n\cdot\Gamma_{n,q}
= -2\lambda < 0.
$$
Thus the rays are pairwise distinct.

\paragraph{Incidence on $X'$.}
Now, we take the strict transforms of $I_n$ on $X'$ under
$\beta:S^{[2]}\dashrightarrow X'$ and denote them by $\widetilde I_n$.

\begin{lemma}\label{lem:section-divisors}
	Strict transform under $\beta$ identifies the prime divisors of
	$S^{[2]}$ with the prime divisors of $X'$ other than the $R_{ij}$.
	In particular, the strict transforms $\widetilde I_n\subset X'$ are
	pairwise distinct $G$-invariant prime divisors, and
	$|m\widetilde I_n|=\{m\widetilde I_n\}$ for every $m\geq1$.
\end{lemma}

\begin{proof}
	Write $R=\sum R_{ij}$ and recall that
	$\operatorname{div}_{X'}(\Theta)=R$.
	Since $p^*B=3R$ and the ramification divisor of $p$ is $2R$,
	Hurwitz gives $K_{X'}-R=p^*(K_X+B/3)$ for compatible canonical
	divisors. The simple normal crossings formula gives discrepancy at
	least $1/3$ for every exceptional divisor over $(X,B/3)$.
	The same bound holds over $(X',-R)$ by the well known formulae from~\cite[Proposition~5.20]{KollarMori1998}.
	Computing the latter discrepancies with
	$K_{X'}=\operatorname{div}_{X'}(\Theta)=R$, we have $K_{X'}-R=0$.
	Hence $a(E;X',-R)=\operatorname{ord}_E(\Theta)$ for every prime divisor $E$ exceptional over $X'$.
	
	Thus $\Theta$ has nonnegative order at every divisorial valuation
	of the common function field. Its order is zero precisely at the
	primes of $X'$ other than the $R_{ij}$. On $S^{[2]}$ it is
	therefore a regular four-form, hence nowhere zero since $S^{[2]}$
	is projective with trivial canonical bundle. Every exceptional
	divisor over $S^{[2]}$ has positive discrepancy since $S^{[2]}$ is smooth, so
	the valuations of order zero are exactly the primes of $S^{[2]}$.
	Thus, for every prime divisor $D\subset S^{[2]}$, the valuation
	$\operatorname{ord}_D$ has order zero on $\Theta$. By the calculation
	on $X'$, its centre on $X'$ is a prime divisor distinct from every
	$R_{ij}$. Conversely, every prime divisor of $X'$ distinct from the $R_{ij}$
	has order zero on $\Theta$, so its centre on $S^{[2]}$ is a prime
	divisor.
	
	Equivariance of $\beta$ gives $G$-invariance of $\widetilde I_n$.
	Since every prime divisor of $S^{[2]}$ has a prime strict transform
	on $X'$ with the same order of vanishing, we have an inclusion
	$$
	H^0\!\left(X',\mathcal O_{X'}(m\widetilde I_n)\right)
	\subseteq
	H^0\!\left(S^{[2]},\mathcal O_{S^{[2]}}(mI_n)\right)
	$$
	inside the common function field. By \eqref{eq:rigid} the right-hand
	side consists of constants, so $|m\widetilde I_n|=\{m\widetilde I_n\}$.
\end{proof}

\section{Descent to \texorpdfstring{$\overline M_{0,7}$}{the stable seven-point space} and proof of the main theorem}
\label{sec:descent}

For each $n\geq1$, let $D_n\subset X$ be the image of
$\widetilde I_n\subset X'$ under $p:X'\to X$, with its reduced
structure. Since $p$ is finite, $D_n$ is a prime divisor.
The prime divisors of $X'$ above $D_n$ form a single $G$-orbit. By
Lemma~\ref{lem:section-divisors}, $G$ preserves $\widetilde I_n$
setwise. Hence $\widetilde I_n$ is the unique prime divisor above
$D_n$.

The same lemma gives $\widetilde I_n\neq R_{ij}$ for every $i,j$.
Thus $p$ is \'{e}tale at the generic point of $\widetilde I_n$
by Lemma~\ref{lem:pair-cubic-cover}, and its ramification index there
is one. Consequently,
\begin{equation}\label{eq:descent-pullback}
	p^*D_n=\widetilde I_n .
\end{equation}
Over a general point of $D_n$, all three points of the fibre lie
on $\widetilde I_n$, so the restriction
$\widetilde I_n\to D_n$ has degree three.
Since the $\widetilde I_n$ are pairwise distinct, so are the $D_n$
by \eqref{eq:descent-pullback}.

\begin{proof}[Proof of Theorem~\ref{thm}]
	Fix $n,m\geq1$ and let $D'\in|mD_n|$. By
	\eqref{eq:descent-pullback} and Lemma~\ref{lem:section-divisors}, we get
	$$
	p^*D'\in|m\widetilde I_n|=\{m\widetilde I_n\}.
	$$
	Pushing forward and using $p_*p^*A=3A$ gives $D'=mD_n$.
	Thus $|mD_n|=\{mD_n\}$ for every $m\geq1$.
	The classes $[D_n]$ span pairwise distinct rays in $N^1(X)_{\mathbb R}$, because
	two classes on the same ray would satisfy $aD_n\sim bD_j$
	for some positive integers $a,b$ and then
	$bD_j\in|aD_n|=\{aD_n\}$ would give 
	$D_j=D_n$.
	
	Let $\sigma_n$ be the section of $\mathcal O_X(D_n)$ cutting out
	$D_n$. It spans the degree-$[D_n]$ piece of the Cox ring.
	A nonzero product of two homogeneous sections of nonzero degrees
	has zero divisor equal to the sum of two nonzero effective divisors.
	It therefore cannot be a scalar multiple of $\sigma_n$, whose zero
	divisor is prime. Every homogeneous generating set of the Cox ring
	must contain a generator in each degree $[D_n]$. This proves that
	$\operatorname{Cox}(X)$ is not finitely generated and hence that
	$X$ is not a Mori Dream Space~\cite[Proposition~2.9]{HuKeel2000}.
	
	To prove nonpolyhedrality, associate to any divisor $D$ on $X$ a
	divisor on $S^{[2]}$ as follows. In $p^*D$, discard the
	components $R_{ij}$ and take the strict transform of every remaining
	prime, with its coefficient, under $\beta^{-1}:X'\dashrightarrow S^{[2]}$.
	For a principal divisor $\operatorname{div}_X(h)$, the result is
	$\operatorname{div}_{S^{[2]}}(p^*h)$, since the corresponding
	primes have the same valuation. As
	$\operatorname{Pic}(X)_{\mathbb R}=N^1(X)_{\mathbb R}$, this defines a linear map
	$$
	\Phi:N^1(X)_{\mathbb R}\longrightarrow N^1(S^{[2]})_{\mathbb R}.
	$$
	The map sends effective divisors to effective divisors or zero, and
	therefore sends pseudoeffective classes to pseudoeffective classes.
	By \eqref{eq:descent-pullback}, $\Phi([D_n])=[I_n]$.
	Fix $n\geq1$. The curves $\Gamma_{n,q}\subset I_n$ from
	Section~\ref{sec:elliptic-incidence} have a common numerical class
	$\gamma_n$ as $q$ varies in $S\setminus C_n$. Define
	$$
	\ell_n:N^1(X)_{\mathbb R}\longrightarrow\mathbb R,\qquad
	\ell_n(\alpha)=\Phi(\alpha)\cdot\gamma_n .
	$$
	Then $\ell_n([D_n])=I_n\cdot\gamma_n=-2$. For every other prime
	divisor $E\subset X$, we have $\ell_n([E])\geq0$. Indeed,
	$p^*E$ has no component $\widetilde I_n$, since
	$p(\widetilde I_n)=D_n$. Hence the effective divisor on $S^{[2]}$
	representing $\Phi([E])$ in the construction above has no
	component $I_n$. As the curves $\Gamma_{n,q}$ cover a dense
	open subset of $I_n$, a general one meets this divisor
	nonnegatively.
	
	Choose an ample divisor $A$ on $X$ and write
	$v_E=[E]/(A^3\cdot E)$ for the normalized class of a prime
	divisor $E$. Since $A$ is ample, the slice of
	$\overline{\operatorname{Eff}}(X)$ cut out by $A^3\cdot\alpha=1$
	is compact.
	It is therefore the convex hull of $v_{D_n}$ and the closed
	convex hull of the $v_E$ for $E\neq D_n$. On the latter set
	$\ell_n$ is nonnegative, while $\ell_n(v_{D_n})<0$.
	Thus $v_{D_n}$ is the unique minimizer of $\ell_n$ on the
	slice, so the $\mathbb R_{\geq0}[D_n]$ are pairwise different extremal rays of
	$\overline{\operatorname{Eff}}(X)$ and consequently $\overline{\operatorname{Eff}}(X)$ is not polyhedral.
\end{proof}

The appendix proves Theorem~\ref{thm} in positive characteristic.

\subsection{A Cantor set on a plane section of the cone}\label{sec:cantor}

Let $F\in N^1(S)_{\mathbb R}$ be the class of a fibre
of the elliptic fibration $f:S\to\mathbb P^1_u$ from
Section~\ref{sec:elliptic-incidence}. The rays
$\mathbb R_{\geq0}[C_n]$ converge to $\mathbb R_{\geq0}F$. Indeed, the degrees
of the $C_n$ against a fixed ample class tend to infinity by the Hodge
index theorem, since their classes are distinct and $C_n^2=-2$. As
$C_n\cdot F=1$, every limit after division by degree has self-intersection
zero and is orthogonal to $F$, hence lies on the ray of $F$.

Recall from Section~\ref{sec:hilbert} that
$Q=\mathbb P^1_x\times\mathbb P^1_y$ and that the triple cover
$S_0\to Q$ is branched along the six lines
$x=1,0,\infty$ and $y=1,0,\infty$.
Write $\varphi:S\to Q$ for the composite $S\to S_0\to Q$.
Let $A_0$ be the class of $\varphi^*\{x=2\}$ and $B_0$ the class
of $\varphi^*\{y=2\}$. These are the fibre classes of the maps
$x\circ\varphi$ and $y\circ\varphi$ from $S$ to $\mathbb P^1$.

To index the six lines, put $c_1=1$, $c_2=0$, $c_3=\infty$.
Let $A_i$ be the unique irreducible curve on $S$ mapping onto
$\{x=c_i\}$, and let $B_j$ be the one mapping onto $\{y=c_j\}$.
Over $(c_i,c_j)\in Q$, the resolution $S\to S_0$ has two
exceptional $(-2)$-curves. Denote the one meeting $A_i$ by
$U_{ij}$ and the one meeting $B_j$ by $V_{ij}$.

The classes $A_0,B_0$ and the eighteen exceptional curves form a basis
of $N^1(S)_{\mathbb R}$. Indeed, $A_0^2=B_0^2=0$ and $A_0\cdot B_0=3$
because $\varphi$ has degree three. Both classes are orthogonal to
the exceptional curves, whose intersection matrix consists of nine
blocks $\bigl(\begin{smallmatrix}-2&1\\1&-2\end{smallmatrix}\bigr)$.
This gives $20$ independent classes, the maximum for a complex K3
surface. The deck group $G$ fixes $A_0,B_0$ because it acts over $Q$.
It also preserves each exceptional pair, and a group of order three
cannot interchange its two curves.

For $\sigma_0=\mathrm{id}$, $\sigma_1=(12)$ and $\sigma_2=(13)$ put
$$
\alpha_j=2A_0+B_0-\sum_{i=1}^3(U_{i,\sigma_j(i)}+2V_{i,\sigma_j(i)}),
\qquad s_j(v)=v+\tfrac13(v\cdot\alpha_j)\alpha_j .
$$
On $V=\langle\alpha_0,\alpha_1,\alpha_2\rangle_{\mathbb R}$ the
intersection matrix is $\left(\begin{smallmatrix}-6&6&6\\6&-6&12\\6&12&-6\end{smallmatrix}\right)$.
Utsumi's fibre components give $3F=\alpha_0+\alpha_1$. The fibre class,
zero section $V_{11}$ and components of the two reducible fibres span
$\alpha_1^\perp$~\cite[\S4]{Utsumi2016}.
Fibrewise inversion fixes this hyperplane pointwise because both component
groups have exponent two and the multiplicity-one components determine
the action on each dual graph. It acts by $-1$ on the Mordell-Weil
group of rank one, hence induces $s_1$. Permutations of the branch values
on the second factor lift to $S$ and give $s_0,s_2$ by conjugation.
Write $\Gamma=\langle s_0,s_1,s_2\rangle$.

The divisor $C\times S+S\times C$ descends to the symmetric
square~\cite[Lemma~7.1.1(iv)]{BrionKumar2005} and pulls back from there
to $I_C$ on $S^{[2]}$, defining a linear map $\mu:[C]\mapsto[I_C]$.
On the normalized graph $X'\xleftarrow{\pi}M\xrightarrow{\nu}S^{[2]}$ of
$\beta$ set $T=\tfrac13p_*\pi_*\nu^*\mu$. Since $\tfrac13p^*p_*$
averages over $G$, Lemma~\ref{lem:section-divisors} gives $\Phi T=\mu$.
For $q\in S$, let $\iota_q:S\setminus\{q\}\to S^{[2]}$ send $a$ to
$\{a,q\}$. Pulling back a line bundle along $\iota_q$ and extending
it over $q$ defines a map $J:N^1(S^{[2]})_{\mathbb R}\to N^1(S)_{\mathbb R}$.
The extension is unique since $S$ is smooth, and its numerical class
is independent of $q$ by algebraic equivalence. For a curve $C\subset S$,
choosing $q\notin C$ gives $\iota_q^*I_C=C|_{S\setminus\{q\}}$,
so $J\mu=\mathrm{id}$.

Now fix an effective divisor $D$ on $S^{[2]}$. Since pairs of distinct
points are dense, we can choose a pair $\{a,q\}$ outside the support
of $D$. For this $q$, the image of $\iota_q$ is not contained in $D$,
so $\iota_q^*D$ is effective, possibly zero. Its closure in $S$
represents $J([D])$. Thus $J$, like $\Phi$ and $T$, preserves
pseudoeffectivity. Since $J\Phi T=J\mu=\mathrm{id}$, we obtain
\begin{equation}\label{eq:cantor-transfer}
	\overline{\Eff}(X)\cap T(V)=T\bigl(\overline{\Eff}(S)\cap V\bigr).
\end{equation}

Put $h=\alpha_0+\alpha_1+\alpha_2$ and $e=(\alpha_1+\alpha_2-\alpha_0)/3$.
In $\mathcal P=\{v\in V:h\cdot v=1\}$,
the region $v^2\geq0$ is a disk for the metric given by minus the
intersection form on $V\cap h^\perp$. Write $\widehat v=v/(h\cdot v)$
for $h\cdot v>0$. Let $K$ be the closure of
$\{\widehat{\gamma F}:\gamma\in\Gamma\}$ on the circle $v^2=0$.

\begin{theorem}\label{thm:cantor-section}
The plane section $\overline{\Eff}(X)\cap T(\mathcal P)$ is the image
under $T$ of the intersection of the half-planes $v\cdot k\geq0$, $k\in K$.
The set $K$ is a Cantor set. For every $k\in K$, the ray spanned by
$T(k)$ lies on the boundary of $\overline{\Eff}(X)$ and is a limit
of rays spanned by extremal prime divisors.
\end{theorem}

\begin{figure}[!t]
	\centering
	\begingroup
	\def\sectionOutline{(0.918559,-0.395285) (1.088662,0.000000) (0.918441,0.395559) (0.912196,0.409767) (0.908238,0.418454)
		(0.903514,0.428558) (0.899197,0.437543) (0.892254,0.451730) (0.886002,0.463680) (0.879304,0.476280)
		(0.873114,0.487516) (0.867727,0.497040) (0.862237,0.506505) (0.856643,0.515909) (0.849023,0.528356)
		(0.840529,0.541766) (0.834119,0.551663) (0.824501,0.565861) (0.819425,0.573187) (0.813983,0.580889)
		(0.768467,0.644379) (0.710989,0.703203) (0.703295,0.710908) (0.697542,0.716545) (0.689978,0.723831)
		(0.683273,0.730163) (0.677477,0.735551) (0.670607,0.741813) (0.645602,0.764264) (0.614591,0.788846)
		(0.602538,0.798092) (0.584652,0.811441) (0.564146,0.825675) (0.548486,0.836218) (0.530087,0.847943)
		(0.519905,0.854224) (0.509649,0.860383) (0.499319,0.866418) (0.488919,0.872329) (0.478448,0.878116)
		(0.467909,0.883777) (0.457303,0.889311) (0.446631,0.894718) (0.435896,0.899997) (0.425098,0.905147)
		(0.414239,0.910168) (0.403321,0.915059) (0.392345,0.919818) (0.381313,0.924446) (0.370227,0.928941)
		(0.359087,0.933304) (0.347897,0.937533) (0.336656,0.941628) (0.325367,0.945588) (0.314032,0.949412)
		(0.302651,0.953101) (0.291228,0.956654) (0.279763,0.960069) (0.268257,0.963347) (0.245672,0.969420)
		(0.227333,0.973817) (0.220088,0.975480) (0.194984,0.980996) (0.166429,0.986053) (0.152798,0.988257)
		(0.144076,0.989567) (0.098969,0.995859) (0.057373,0.998353) (0.042554,0.999103) (0.031777,0.999495)
		(0.020816,0.999783) (0.009852,0.999951) (0.001748,0.999998) (-0.012786,0.999918) (-0.027233,0.999629)
		(-0.148454,0.995859) (-0.262022,0.965062) (-0.272007,0.962295) (-0.283535,0.958962) (-0.295421,0.955368)
		(-0.304797,0.952417) (-0.322567,0.946681) (-0.347804,0.937567) (-0.356008,0.934503) (-0.367585,0.929990)
		(-0.384352,0.923187) (-0.395086,0.918644) (-0.405766,0.913977) (-0.416391,0.909185) (-0.426960,0.904270)
		(-0.439540,0.898226) (-0.448950,0.893566) (-0.457270,0.889328) (-0.466569,0.884521) (-0.474355,0.880334)
		(-0.486914,0.873450) (-0.512312,0.859173) (-0.535774,0.844361) (-0.545873,0.837876) (-0.555771,0.831336)
		(-0.564178,0.825653) (-0.571909,0.820317) (-0.579143,0.815226) (-0.587234,0.809442) (-0.595325,0.803486)
		(-0.604097,0.796911) (-0.612372,0.790569) (-1.632993,0.000000) (-0.612372,-0.790569) (-0.601826,-0.798627)
		(-0.587234,-0.809442) (-0.576643,-0.816999) (-0.564178,-0.825653) (-0.552733,-0.833360) (-0.542802,-0.839861)
		(-0.534732,-0.845022) (-0.512312,-0.859173) (-0.488604,-0.872506) (-0.478757,-0.877948) (-0.466569,-0.884521)
		(-0.457669,-0.889123) (-0.448950,-0.893566) (-0.441605,-0.897210) (-0.426960,-0.904270) (-0.416391,-0.909185)
		(-0.405766,-0.913977) (-0.395086,-0.918644) (-0.384352,-0.923187) (-0.370579,-0.928806) (-0.356008,-0.934503)
		(-0.342224,-0.939618) (-0.322567,-0.946681) (-0.302817,-0.953049) (-0.294694,-0.955592) (-0.283984,-0.958829)
		(-0.275308,-0.961361) (-0.266877,-0.963731) (-0.148454,-0.995859) (-0.027233,-0.999629) (-0.013877,-0.999904)
		(-0.000065,-1.000000) (0.009852,-0.999951) (0.020816,-0.999783) (0.031777,-0.999495) (0.042554,-0.999103)
		(0.057373,-0.998353) (0.098969,-0.995859) (0.140096,-0.990138) (0.152798,-0.988257) (0.168277,-0.985742)
		(0.194984,-0.980996) (0.220088,-0.975480) (0.227333,-0.973817) (0.245672,-0.969420) (0.268257,-0.963347)
		(0.279763,-0.960069) (0.291228,-0.956654) (0.302651,-0.953101) (0.314032,-0.949412) (0.325367,-0.945588)
		(0.336656,-0.941628) (0.347897,-0.937533) (0.359087,-0.933304) (0.370227,-0.928941) (0.381313,-0.924446)
		(0.392345,-0.919818) (0.403321,-0.915059) (0.414239,-0.910168) (0.425098,-0.905147) (0.435896,-0.899997)
		(0.446631,-0.894718) (0.457303,-0.889311) (0.467909,-0.883777) (0.478448,-0.878116) (0.488919,-0.872329)
		(0.499319,-0.866418) (0.509649,-0.860383) (0.519905,-0.854224) (0.530087,-0.847943) (0.548486,-0.836218)
		(0.564146,-0.825675) (0.584652,-0.811441) (0.602538,-0.798092) (0.614591,-0.788846) (0.645602,-0.764264)
		(0.669686,-0.742644) (0.677477,-0.735551) (0.689978,-0.723831) (0.697542,-0.716545) (0.703295,-0.710908)
		(0.713311,-0.700848) (0.768467,-0.644379) (0.813983,-0.580889) (0.819425,-0.573187) (0.824501,-0.565861)
		(0.834119,-0.551663) (0.840529,-0.541766) (0.849023,-0.528356) (0.856643,-0.515909) (0.862237,-0.506505)
		(0.867727,-0.497040) (0.874147,-0.485665) (0.879304,-0.476280) (0.882253,-0.470776) (0.892254,-0.451730)
		(0.901374,-0.433046) (0.908238,-0.418454) (0.912196,-0.409767)}
	\definecolor{outlineColor}{RGB}{28,57,66}
	\definecolor{sectionColor}{RGB}{207,224,223}
	\definecolor{diskColor}{RGB}{244,248,247}
	\definecolor{circleColor}{RGB}{94,111,115}
	\definecolor{arcColor}{RGB}{148,65,28}
	\resizebox{.96\linewidth}{!}{%
		\begin{tikzpicture}[x=1cm,y=1cm,
			every node/.style={inner sep=1.5pt,text=black,font=\small}]
			\begin{scope}[x=2.04cm,y=2.04cm]
				\path[fill=sectionColor] plot coordinates {\sectionOutline};
				\fill[diskColor] (0,0) circle[radius=1];
				\draw[circleColor,line width=.48pt,dash pattern=on 2pt off 1.55pt]
				(0,0) circle[radius=1];
				\draw[outlineColor,line width=.8pt,line join=round]
				plot coordinates {\sectionOutline} to cycle;
				\foreach \P in {
					(-0.599772,-0.800171),
					(-0.563771,-0.825931),
					(-0.408248,-0.912871),
					(-0.288761,-0.957401),
					(0.019440,-0.999811),
					(0.408248,-0.912871),
					(0.861858,-0.507151),
					(0.906112,-0.423038),
					(0.906112,0.423038),
					(0.861858,0.507151),
					(0.408248,0.912871),
					(0.019440,0.999811),
					(-0.288761,0.957401),
					(-0.408248,0.912871),
					(-0.563771,0.825931),
					(-0.599772,0.800171)}
				{\path[fill=arcColor] \P circle[radius=1.05pt];}
				\draw[arcColor,line width=.54pt]
				(80:1) to (80:.85) arc[start angle=80,end angle=115,radius=.85] to (115:1);
			\end{scope}
			\begin{scope}[shift={(4.04,0)},x=.21cm,y=1cm]
				\node at (17.5,1.11) {Polar angle on the marked arc};
				\draw[circleColor,line width=.46pt] (0,.42) to (35,.42);
				\foreach \X/\A in {0/80,5/85,10/90,15/95,20/100,25/105,30/110,35/115}{
					\draw[circleColor,line width=.38pt] (\X,.42) to (\X,.49);
					\node[anchor=south,font=\scriptsize] at (\X,.50) {$\A^\circ$};
				}
				\def\angleRow#1#2#3{%
					\node[anchor=east] at (-.68,#2) {$d=#1$};
					\draw[circleColor!55,line width=.32pt] (0,#2) to (35,#2);
					\foreach \firstAngle/\lastAngle in {#3}
					{\fill[outlineColor] (\firstAngle,{#2-.054}) rectangle (\lastAngle,{#2+.054});}%
				}
				\angleRow{2}{.06}{0.000000/11.784923,25.172504/35.000000}
				\angleRow{4}{-.40}{0.000000/2.084207,6.564864/11.784923,25.172504/27.902616,29.728845/35.000000}
				\angleRow{7}{-.86}{0.000000/1.395276,1.647298/2.084207,6.564864/6.745036,6.840665/7.320207,
					7.802084/9.638871,9.843487/10.081308,10.695149/10.822129,10.907593/11.367152,
					11.509263/11.675470,11.712152/11.784923,25.172504/25.216312,25.238323/25.337465,
					25.421473/25.688442,25.737314/25.809488,26.151143/26.280374,26.390201/27.323214,
					27.553287/27.776625,27.820510/27.902616,29.728845/29.887356,29.977396/30.504268,
					31.205383/35.000000}
			\end{scope}
		\end{tikzpicture}%
	}%
	\endgroup
	\caption{The dashed circle $v^2=0$ bounds a disk in $\mathcal P$.
		Let $H_d$ be the convex hull of this disk and the points obtained
		from $e=(\alpha_1+\alpha_2-\alpha_0)/3$ by at most $d$ reflections from
		$\{s_0,s_1,s_2\}$.	The left and right tips are $\widehat e$ and $\widehat{s_0(e)}$, respectively. The solid line shows the boundary of $H_9$, and the dots mark points of the Cantor set $K$.
		The right panel displays the marked arc from $80^\circ$ to $115^\circ$
		by polar angle. In row $d$, the dark intervals show where the arc
		meets the boundary of $H_d$. A point of the arc belongs to $K$
		exactly when its angle is dark for every $d\geq 0$.}
	\label{fig:cantor}
\end{figure}
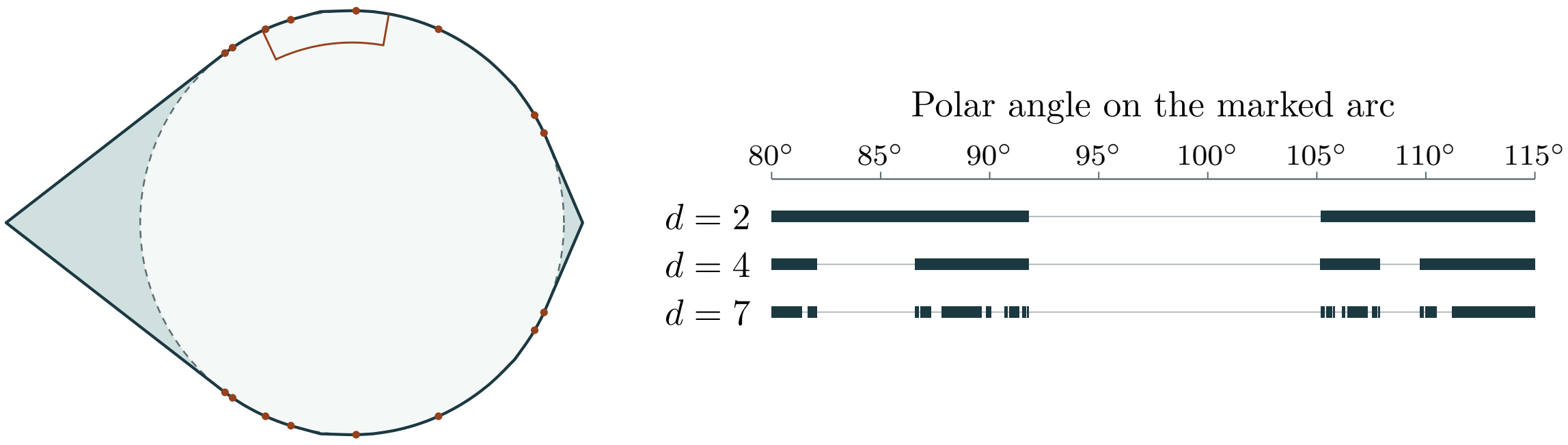

\begin{proof}
In the Klein model of the hyperbolic plane, geodesics are straight
segments in the disk.
The matrix above shows that the lines $v\cdot\alpha_j=0$ and $v\cdot e=0$
bound a quadrilateral with two vertices on the circle and two right angles.
Gluing two copies along the reflection sides gives a surface with
two cusps and one geodesic boundary component. The deck transformations of
its universal cover are the even products of reflections. Since $\widehat F$
represents a cusp, its orbit closure $K$ is the limit
set~\cite[Theorem~5.3.7]{Beardon1983}. The boundary lines of the
hyperbolic convex hull of $K$ are the lifts $v\cdot\gamma e=0$
of the geodesic boundary~\cite[proof of Lemma~4.4, p.~634]{AlessandriniEtAl2011}.
The limit set is perfect. Its boundary in the circle is nonempty and
invariant, hence equals $K$ by minimality. Thus $K$ is a Cantor set.

Every $k\in K$ is a limit of effective nef fibre classes. Hence
$v\cdot k\geq0$ for every pseudoeffective $v$. Conversely, the disk
centre $h/30$ is inside the quadrilateral and hence inside the convex
hull of $K$. At the endpoints of each complementary arc the two tangents
meet at $\widehat{\gamma e}$ for the corresponding boundary line.
This point is pseudoeffective since
$$ e=2A_1+B_1+2U_{11}+2V_{11}+U_{12}+U_{13}. $$
Each complementary arc contributes the two tangent segments from its
endpoints to $\widehat{\gamma e}$. Thus the intersection of the tangent
half-planes is the closed convex hull of $K$ and these points.
Equation~\eqref{eq:cantor-transfer} gives the plane section in Figure~\ref{fig:cantor}.
Its boundary points $T(k)$ also lie on the boundary of the full cone.

Every smooth rational curve $C\subset S$ is $G$-invariant by $C^2=-2$,
so $I_C$ gives an extremal prime divisor $D_C$ on $X$ as in
Section~\ref{sec:descent}. Unless $C$ is one of the finitely many curves
fixed pointwise, its fixed points form a finite set. A fixed reduced
pair in $I_C$ consists of such a point and a point of $S^G$, of dimension
at most one. A nonreduced fixed pair is supported at a fixed point of
$C$, where its tangent direction varies in a projective line.
Thus $\dim I_C^G\leq1$. Each $R_{ij}$ has discrepancy one over $S^{[2]}$
because $\operatorname{ord}_{R_{ij}}(\Theta)=1$, so its centre is a surface
fixed pointwise by $G$. Hence $I_C$ contains none of these centres and
$\pi_*\nu^*I_C$ is its strict transform, of degree three over $D_C$.
Thus $T([C])=[D_C]$ for all but finitely many $C$. Applying the
automorphisms inducing $\gamma$ to the $C_n$ gives extremal rays tending
to the ray of $T(\gamma F)$, and the points $\widehat{\gamma F}$ are
dense in $K$.
\end{proof}

\section{Appendix: positive characteristic}\label{sec:positive}

	The proof of Theorem~\ref{thm} also works over an algebraically closed
	field of characteristic $p\geq5$.  $S_0$ still
	has nine $A_2$ singularities, so its resolution $S$ is again a K3
	surface and $S^{[2]}$ has trivial canonical bundle~\cite[\S4.3, Corollary~1]{KumarThomsen2001}.
	Identity \eqref{eq:flex-identity} holds over $\mathbb Z$.
	With $x_0=0,\ x_1=1$, we have
	$$
	q(P)=x_2(x_2-1)(x_4-x_3)(x_5-x_3)(x_5-x_4)\neq0
	$$
	on $M_{0,7}$. Let $m$ be the projection of $\ell(x,\tau)$ from $o$.
	Then $q|_m$ cuts out $m\cap S_0$ away from $t=0$, where $z=t\tau\neq0$.
	Its span and the inverse formulas of Sections~\ref{sec:field}
	and~\ref{sec:hilbert} recover the marked curve and its cube root.
	
	The pencil \eqref{eq:weier} has smooth generic fibre of genus one.
	Adjunction on a resolution of the pencil excludes a horizontal
	exceptional curve, so the pencil has no base point on $S$.
	Over the $k$-line, $P=(1,-1)$ specializes to points of orders
	four and three at the smooth fibres $k=-1/4$ and $k=-3/4$,
	respectively. Injectivity of good reduction on prime-to-$p$
	torsion proves that $P$ has infinite order~\cite[Proposition~3.1]{Silverman2009}.
	Finally, the cyclic extension of degree three is tame at every
	divisorial valuation, so the discrepancy calculation in
	Lemma~\ref{lem:section-divisors} and the descent of the incidence
	divisors apply as well~\cite[Def.~4.1 and proof of Lemma~4.3]{Stigant2023}.
	
	In characteristic three, $\tau^3=\delta$ defines a purely inseparable
	cover and $p^*\theta=0$. In characteristic two, the cubic cover
	remains tame, but the pencil used in Section~\ref{sec:elliptic-incidence}
	has geometrically rational generic fibre. Thus, we have to adapt our proof in these cases.

\subsection{Characteristic two}\label{sec:char-two}

Let $\Bbbk$ be algebraically closed of characteristic two and put
$H=S^{[2]}$. The cubic cover stays tame. The arguments from Sections~\ref{sec:field} and~\ref{sec:hilbert} again identify $X'$
birationally with $H$. Indeed, project the line $\ell(x,\tau)$ from $o$
to a line $m$ in the space of directions. Since $q(P)\neq0$, the
length-two subscheme $m\cap S_0$ cut out by $q|_m$ lies in the smooth
open subset $z\neq0$ of $S_0$. Its linear span recovers $m$.
 We construct the required $G$-invariant divisors on $H$ using curves of genus two.

On $d=1$, write $x=b/d,\ y=c/d$, so that
$z^3=x(x+1)y(y+1)$. Let $A,B\in S_0$ be the unique points above $(0,0),(\infty,\infty)\in Q$, respectively. Project $S_0$ from the line joining $A$ and $B$. The projection
$[a:b:c:d:z]\mapsto[b:c:z]$ extends to
$\pi:S\to\mathbb P^2$: at each base point its three coordinates
generate the maximal ideal, and blowing up resolves the $A_2$ point.
The two exceptional curves over $A$ map isomorphically to lines.
Denote them by $A_x,A_y$, indicating if they meet the strict transform
of $x=0$ or $y=0$, and label $B_x,B_y$ similarly at infinity.
Above $[1:t:s]$, substituting $y=tx,\ z=sx$ gives
$$
t^2x^2+(t^2+t+s^3)x+t=0 .
$$
Thus $\pi$ has separable degree two. Its contracted curves are the
four branch rulings through $A$ or $B$ and the fourteen exceptional
curves over the other seven crossings. Indeed, the displayed polynomial
vanishes identically only for $s=t=0$; interchanging $x,y$ and
examining $d=0$ gives the other three rulings.

Let $C$ be the inverse image of the generic line $z=\alpha x+\beta y$,
over $K=\Bbbk(\alpha,\beta)$. This curve is geometrically smooth.
Away from $z=0$, a singular point would satisfy
$$
y(y+1)=\alpha z^2,\qquad x(x+1)=\beta z^2.
$$
Multiply the first equation by $x$ and the second by $y$. Adding gives
$xy(x+y)=(\alpha x+\beta y)z^2=z^3$. The surface equation gives
$z^3=xy(x+1)(y+1)$. Since $z\neq0$, we can cancel $xy$.
Expanding $(x+1)(y+1)=x+y$ in characteristic two gives $xy=1$. These points form at most a curve on $S_0$, and each determines $(\alpha,\beta)$.
At $z=0$ the generic section is smooth except for ordinary nodes
at $A,B$, which the resolution separates. Bertini irreducibility
and adjunction now give a smooth genus-two curve, since
$h=\pi^*\mathcal O(1)$ has $h^2=2$.

The points $p_x=C\cap A_x$ and $p_y=C\cap A_y$ are $K$-rational.
We now claim that $p_x-p_y$ has infinite order. 
For this, let $V_\pi$ be the rational span of the eighteen contracted curves. The universal family of inverse images of lines is a $\mathbb P^1$-bundle over $S$. A class trivial on the generic curve is represented there by vertical divisors. Modulo
the pullbacks of contracted curves, these are rational combinations
of pullbacks from the dual plane: only finitely many lines split into
two components dominating the line. Each such component maps
isomorphically to the line and has $h$-degree one and square $-2$;
the Hodge index theorem gives finitely many such classes, each with
a unique irreducible representative. The projective bundle formula
gives that the kernel of $\operatorname{Pic}(S)_{\mathbb Q}\to\operatorname{Pic}(C)_{\mathbb Q}$ is contained in $V_\pi$.

Now let $F_x,F_y$ be the two fibre classes pulled back from $Q$.
The class $L_\pi=A_x-A_y+B_x-B_y+F_x-F_y$ is orthogonal to all eighteen
contracted curves, whereas $L_\pi\cdot(A_x-A_y)=-6$. Thus $A_x-A_y$
lies outside $V_\pi$, proving the claim that $p_x-p_y$ has infinite order. 

Two general points of $S$ determine a line through their images
under $\pi$. Hence $H$ is birational over the dual plane to the
relative $\operatorname{Pic}^2$ of this family~\cite[Theorem~5.1]{Milne1986Jacobians}. Translation by $p_x-p_y$
commutes with $G$. Its iterates send $I_{A_x}$ to distinct prime
divisors: on the generic fibre they translate the theta divisor
$p_x+C$, whose translation stabilizer is trivial~\cite[Theorem~6.6]{Milne1986Jacobians}.
Since $H$ is smooth with trivial canonical bundle~\cite[Thms.~7.4.1 and~7.4.6]{BrionKumar2005}, these birational
translations are isomorphisms in codimension one, by the volume-form
argument of Lemma~\ref{lem:section-divisors}. They therefore transport
both $I_{A_x}$ and its separating intersection functional.
The resulting divisors are $G$-invariant, so the tame descent of
Section~\ref{sec:descent} gives infinitely many rigid prime divisors
on $X$ spanning distinct extremal rays.

\begin{remark}
	The Jacobian translations above also give divisors on $X$ over $\mathbb C$.
	For the same projection $[a:b:c:d:z]\mapsto[b:c:z]$, the inverse image
	of a general line is smooth by Bertini. In the coordinates $y=tx$, $z=sx$
	the equation is $t^2x^2+(s^3-t^2-t)x+t=0$. Consequently the projection
	has the same eighteen contracted curves, and the intersection computation
	with $A_x-A_y$ proves that $p_x-p_y$ has infinite order as above.
	These translations give infinitely many divisors distinct from the
	surface incidence divisors of Section~\ref{sec:elliptic-incidence}.
	Indeed, on $\operatorname{Pic}^2(C)$, incidence with a section curve $D$
	of positive $h$-degree gives one translate $p+C$ for each of the
	$h\cdot D$ distinct geometric points $p\in C\cap D$. These translates
	are distinct by the triviality of the theta stabilizer. A curve of
	$h$-degree zero misses the generic member of the family. Every Jacobian
	translate under consideration restricts to a single theta translate,
	so equality requires $h\cdot D=1$. The Hodge index theorem gives only
	finitely many classes with $D^2=-2$ and $h\cdot D=1$, each with at most
	one irreducible representative. The prime-divisor correspondence and
	finite descent preserve distinctness.
\end{remark}

\subsection{Characteristic three}\label{sec:char-three}

Let $\Bbbk$ be algebraically closed of characteristic three and put
$H=S^{[2]}$. The cover $p:X'\to X$ given by $v^3=s_B$ is now
purely inseparable of degree three, since $s_B$ has a simple zero
along each $B_{ij}$. We keep the elliptic sections of
Section~\ref{sec:elliptic-incidence} and describe the change in
their descent.

In the function field of $S_0$ put $x=1/r_1,\ y=1/r_2,\ w=Z/(r_1r_2)$.
In these coordinates the surface equation is $w^3=(x^3-x)(y^3-y)$. This is the model
of Shimada-Zhang: its resolution $S$ is the supersingular K3
surface of Artin invariant one, and $S_0$ has ten $A_2$ points~\cite[Examples~1.7 and~7.8]{ShimadaZhang2007}. Nine of these singularities lie above the intersections of the six
fibres on $Q$ described in Section~\ref{sec:hilbert}. The tenth is
$v_*=(0,0,1)$ in the coordinates $(r_1,r_2,Z)$.
Write $E_1,E_2$ for its two exceptional curves.

We first determine the images of $I_{E_1},I_{E_2}$ under the
birational map $\beta:H\dashrightarrow X'$. Take a primitive fifth root $\epsilon$
and the six finite marks $(0,1,\epsilon,\epsilon^2,\epsilon^3,\epsilon^4)$.
Writing $c=1-\epsilon^2-\epsilon^3$, substitution gives
$$
\delta=c=(-c)^3,\qquad V=c(1,-1,-1,1,-1),\qquad
q(sP+tV)=s^2q(P)+st\,c(1+\epsilon^3).
$$
With $\tau=-c$, the point $[V:\tau]$ projects to $v_*$.
The other root of $q$ is distinct and gives a smooth point of
$S_0$, since $q(P)\neq0$ and $v_*$ is its only singular point
with $z\neq0$. Thus for a general matching line the projected line $m$ meets
$S_0$ in two distinct points. If one is $v_*$, the other point
$u\in S_0$ determines $m=\langle v_*,u\rangle$ and the line above it
in $W$ by $h=g|_m/q|_m$. Such lines therefore form a family of
dimension at most two. The inverse formulas of
Section~\ref{sec:field} therefore give the birational map
$\beta:H\dashrightarrow X'$ as before. For a general pair in
either $I_{E_i}$, the line in $W$ depends only on $v_*$ and the
image in $S_0$ of the other point of the pair. The displayed configuration shows that the matching inverse is defined generically on this family. Hence $\beta$
contracts both $I_{E_i}$.

To show that $I_{E_1}$ and $I_{E_2}$ are the only prime divisors
contracted by $\beta$, we first compute $H^0(X',\omega_{X'})$.
Adjunction gives $\omega_{X'}\simeq p^*\mathcal O_X(-K_X)$.
A general point of $B_{ij}$ represents two copies of $\mathbb P^1$
meeting at a node, one carrying $p_i,p_j$ and the other carrying the
remaining five marks. Keep both components and all marks fixed, and
glue a fixed point of the first component to a variable point of the
second. The closure of this family is a curve $C\subset B_{ij}$
through the given point, with $K_X\cdot C=1$~\cite[\S2.2, before
Lemma~2.4, and Corollary~2.7, Table~1]{Moon2017Seven}.
An effective divisor $D\sim-K_X$ would satisfy $D\cdot C=-1$ and
hence contain $C$. Since such a curve passes through a general point
of each $B_{ij}$, the divisor $D-B$ would be effective. But
$D-B\sim_{\mathbb Q}-2B/3$ has negative intersection with the cube
of an ample divisor. Thus $H^0(X,\mathcal O_X(-K_X))=0$.
Likewise, $K_X\sim_{\mathbb Q}-B/3$ gives
$H^0(X,\mathcal O_X(K_X))=0$. By the projection formula,
$$
p_*\omega_{X'}\simeq
\mathcal O_X(-K_X)\oplus\mathcal O_X\oplus\mathcal O_X(K_X).
$$
Only the middle summand has global sections. Its constant section
$1$ corresponds to the tautological section $v$ of
$p^*\mathcal O_X(-K_X)$, so $H^0(X',\omega_{X'})=\Bbbk v$.

Let $\nu:N\to X'$ be the normalization. A volume form on the
smooth fourfold $H$ extends to $N$; its order is positive
precisely at primes exceptional over $H$. Finite duality identifies
$\nu_*\omega_N$ with the conductor times $\omega_{X'}$~\cite[Tag~0AX4]{StacksProject}. The form thus maps to a nonzero
multiple of $v$, so the conductor is the unit ideal off $v=0$.
Along $v=0$ the hypersurface $X'$ is regular in codimension one
because $B$ is reduced. Thus $X'$ is normal by Serre's criterion,
and the divisor of the form on $X'$ is precisely $\sum R_{ij}$.
Thus $\beta^{-1}$ contracts exactly the twenty-one $R_{ij}$.

Purely inseparable pullback identifies
$\operatorname{Cl}(X')_{\mathbb Q}$ with
$\operatorname{Pic}(X)_{\mathbb Q}$, of rank $42$.
On the other side $\operatorname{Pic}(H)$ has rank
$1+\rho(S)=23$, by Hilbert-Chow. Comparing the divisor
localization sequences on a common isomorphic open subset gives
$21+23-42=2$ primes contracted by $\beta$: the relations
among omitted divisors come from the same units on this open subset.
Thus the missing primes are exactly $I_{E_1}$ and $I_{E_2}$.

Finally, putting $t=-u$ in \eqref{eq:weier} gives Sengupta's
Fibration~5. His section $(1,1)=-P$ has infinite order and meets
the identity component at $t=1$~\cite[Fibration~5; preprint \S7, p.~11]{Sengupta2017}.
Here $v_*$ maps to the node $(-1,0)$ of the Weierstrass fibre
$Y_e^2=X_e(X_e+1)^2$; the coordinate change is regular in both
directions at this point. Consequently every $C_n$ misses both
$E_i$, since points with nonsingular reduction form a subgroup~\cite[Chapter~VII, Proposition~2.1]{Silverman2009}.

Let $D_n$ be the reduced image on $X$ of the prime transform of
$I_n$. Since $p$ is radicial, $p^*D_n=e_n\widetilde I_n$ for
some $e_n>0$. Pull a divisor on $X$ back to $X'$, discard the
$R_{ij}$, take strict transform on $H$, and intersect with
$\Gamma_{n,q}$ for general $q$. For a principal divisor the only
omitted terms on $H$ are multiples of $I_{E_1},I_{E_2}$, and
both pair to zero because $C_n\cdot E_i=0$. We therefore obtain
a numerical functional on $X$. It takes the value $-2e_n$ on
$D_n$ and is nonnegative on every other prime. The arguments of Section~\ref{sec:descent}
now give rigidity, distinct extremal rays, and infinitely many Cox
generator degrees.

\bibliographystyle{alpha}
\bibliography{M07_references}
\end{document}